\documentclass[a4paper,leqno]{article}

\usepackage{amsmath, amssymb}
\input graphicx
\date{}
\title{A Gauss-Bonnet Type Formula on Riemann-Finsler surfaces with non-constant indicatrix volume}
\author{J. Itoh, S.V. Sabau, H. Shimada}

\begin{document}

\maketitle
\begin{center}
{\it 
In the memory of Professor Makoto Matsumoto
}
\end{center}

\begin{abstract}
We prove a Gauss-Bonnet type formula for Riemann-Finsler surfaces of non-constant indicatrix volume and with regular piecewise $C^\infty$ boundary. We give a Hadamard type theorem for $N$-parallels of a Landsberg surface.
\end{abstract}

\bigskip


\bigskip

2000 {\it Mathematics Subjects Classification. } Primary 53C60, Secondary 53C20.

\section{Introduction}

$\quad$ A major topic in Riemannian geometry is the study of the relation between the curvature of the Riemannian metric and the topology of the manifold. This is mainly achieved through the well-known Gauss-Bonnet-Chern theorem. The theorem and its consequences are especially interesting in the case of 
Riemannian surfaces (see \cite{SST2003} for a comprehensive exposition).

The Gauss-Bonnet theorem was extended for the first time by D. Bao and S. S. Chern to the case of 
boundaryless Finsler manifolds of Landsberg type and Finsler manifolds of constant volume (see for details \cite{BC1996}). In the case of Landsberg surfaces the Gauss-Bonnet-Chern theorem is stated in a particular form that can be regarded as a direct generalization to the Finslerian case of the Riemannian classical 
result. In \cite{SS2007} we have extended the Gauss-Bonnet-Chern theorem for boundaryless Landsberg surfaces to the case of Landsberg surfaces with smooth boundary.

The reason to restrict the considerations to Landsberg surfaces is that on these surfaces the Riemannian volume of the indicatrix is constant and therefore the Euler-Poincare characteristic of the manifold can be related to the curvature  in a similar way to the Riemannian case. However, the Landsberg structures include the Berwald ones, which, at least in the case of surfaces, are known to be locally Minkowski in the flat case, or Riemannian otherwise.

Recently, there are many suspicions about the existence of regular Landsberg structures that are not Berwald (\cite{Sz2008a}, \cite{Ma2008}, \cite{Sz2008b}), but the existence of such structures is still an open problem.

However, the Finsler structures more general than the Lansdberg ones can have very interesting geometrical properties, and a Gauss-Bonnet type formula might be a useful tool in the study of their geometry.

In \cite{Sh1996} are proved some Gauss-Bonnet type formulas for $2n$ dimensional Riemann-Finsler manifolds whose indicatrix volume is constant. The paper also contains interesting information on different attempts to extend the Gauss-Bonnet theorem to the Finslerian setting. See also \cite{BS1994} for several discussions on the constancy of the indicatrix volume.

On the other hand, M. Matsumoto studies a Gauss-Bonnet formula for bounded regions on a Finsler surface, but he uses a completely different approach than ours 
(\cite{M1984}).  Matsumoto's normals and curvatures have different geometrical meanings than the ones in the present paper. For his setting S. S. Chern's transgression method used by us can not be employed.  

 In the present paper we are concerned with the following question:

{\it Question. Does a Gauss-Bonnet type formula hold in the case of 
Riemann-Finsler domains with regular piecewise $C^\infty$ 
boundary?}

The lack of angles is a sort of peculiarity of the traditional Finsler geometry. We will show in the present study that the so-called {\it Landsberg angles} can be very useful in the study of the geometry near a ``corner'' of a regular piecewise $C^\infty$ curve.

The paper is organized as follows. We recall some basic facts on the geometry of Riemann--Finsler manifolds in \S2. We discuss here the Landsberg angles defined as the Riemannian length of the indicatrix curve arc defined by the tips of two unit vectors. In \S3 we treat the normal lift of a curve to the indicatrix 
$\Sigma$ which is different from the canonical lift of a curve usually used (see for  example \cite{BCS2000}, p. 112). We are led in this way to the notion of 
$N$-parallels and $N$-parallel curvature of a curve $\gamma$ on the surface $M$. The difference with the Finslerian geodesics is also discussed. An existence and unicity theorem for $N$-parallels is given in the Appendix. 

Theorem 4.2. proved in \S4 gives a partial answer to the question above. We give here a topological lemma that allows us to relate the Euler characteristic of $M$ with the Chern connection 1-form in the case when the indicatrix length is not constant, i.e. in a more general case than the Landsberg structures. 

The regular piecewise $C^\infty$ boundary case is discussed in \S5 where we construct a variation curve near the given boundary. The Gauss lemma for Riemann-Finsler manifolds is the one that makes all the machinery working. Here is where we prove the Theorem 5.1 which gives the final affirmative answer to the question above. 

We finally show how the Gauss-Bonnet theorem controls the behavior of $N$-parallels by proving  a Hadamard type theorem in \S6 for Landsberg surfaces.
\newpage


\section{The geometry of Riemann--Finsler surfaces}

\quad This chapter follows closely \cite{BCS2000}, Chapter 4.

A {\bf Finsler norm}, or metric, on a real, smooth, $2$-dimensional manifold
$M$ is a function $F:TM\to \left[0,\infty \right)$ that is positive and
smooth away from the zero section, has the {\bf homogeneity property}
$F(x,\lambda v)=\lambda F(x,v)$, for all $\lambda > 0$ and all 
$v\in T_xM$, having also the {\bf strong convexity property} that the
Hessian matrix
\begin{equation*}
g_{ij}(x,y)=\frac{1}{2}\frac{\partial^2F^2(x,y)}{\partial y^i \partial y^j}
\end{equation*}
is positive definite at every point of $\widetilde{TM}=TM\setminus \{0\}$.

This implies that the Finslerian unit sphere, or the {\bf indicatrix}
\begin{equation}\tag{2.1}
\Sigma_x:=\{x\in T_xM \ |\  F(x,v)=1\}\subset T_xM
\end{equation}
at $x\in M$ be a smooth, closed, strictly convex hypersurface in
$T_xM$. In addition, if $F(x,-v)=F(x,v)$, then $F$ is said to be
reversible, or absolutely homogeneous (see \cite{M1986}, \cite{BCS2000} or \cite{Sh2001} for the basics of Riemann-Finsler manifolds).
\bigskip

{\bf Remark.}
We gave here the definition for the case of a surface $M$ because in the present paper we deal only with surfaces, but the above definition can be easily extended to the arbitrary dimensional case. 
\bigskip

A smooth 2-dimensional manifold endowed with a Finsler norm is called a {\bf Finsler structure on the surface $M$}, or simply a {\bf Finsler surface}.
\bigskip

In other words, a Finsler surface
 is a pair $(M,F)$ where $F:TM\to [0, \infty)$ is $\mathbf 
C^\infty$ on $\widetilde{TM}:=TM\backslash \{0\}$ and whose restriction to 
each tangent plane $T_xM$ 
is a Minkowski norm (see \cite{SS2007} for a detailed discussion).

A Finsler structure $(M,F)$ on a surface $M$ is also equivalent to a smooth
hypersurface (i.e. 3-dimensional submanifold) $\Sigma\subset TM$ for which the canonical projection
$\pi:\Sigma\to M$ is a surjective submersion and having the property
that for each $x\in M$, the $\pi$-fiber $\Sigma_x=\pi^{-1}(x)$ is
a strictly convex smooth curve including the origin $O_x\in T_xM$. 

Recall that in order to study the geometry of the surface $(M,F)$ one considers the pull-back bundle $\pi^*TM$ with base manifold $\Sigma$ and fibres 
$(T_xM)\vert_u$, where $u\in\Sigma$ such that $\pi(u) = x$ (see \cite{BCS2000}, Chapter 2). In general this is not a principal bundle.
\bigskip

Let us remark that if we denote the projection by $p: TM \longrightarrow M$, then one can start with the pull-back bundle $p^*TM$ constructed over the slit tangent bundle $\widetilde{TM}$. This is also a vector bundle whose fiber over a typical point $u = (x,y)\in \widetilde{TM}$ is a copy of $T_xM$, where $p(x,y) = x \in M$.\\

However, since the majority of our geometrical objects are sections of the pull-back bundle $\pi^*TM$ with base manifold $\Sigma$, we prefer to use this one instead of $p^*TM$ over $\widetilde{TM}$.\\

We point out that we are in fact using the same theory as in \cite{BCS2000}, but we have switched the notation for $p: TM \longrightarrow M$ with $\pi: \Sigma \longrightarrow M$.\\

It is also known (\cite{BCS2000}, p. 30) that the vector bundle $\pi^*TM$ has a distinguished global section $l:=\frac{y^i}{F(y)}\frac{\partial}{\partial x^i}$.\\

Using this section, one can construct a positively oriented $g$-orthonormal 
frame $\{e_1,e_2\}$ for $\pi^*TM$, where $g=g_{ij}(x,y)dx^i\otimes dx^j$ is 
the induced Riemannian metric on the fibers of $\pi^*TM$. The frame $\{u;e_1,e_2\}$ for any $u\in\Sigma$ is a globally defined $g$-orthonormal frame field 
for 
$\pi^*TM$ called the {\bf Berwald frame}. 

Locally, we have,
\begin{equation*}
\begin{split}
e_1&:= \cfrac{1}{\sqrt g}\biggl(\cfrac{\partial F}{\partial 
y^2}\cfrac{\partial}{\partial x^1}-\cfrac{\partial F}{\partial 
y^1}\cfrac{\partial}{\partial x^2}\biggr)=m^1\cfrac{\partial}{\partial 
x^1}+m^2\cfrac{\partial}{\partial x^2}\quad ,\\
e_2&:= \cfrac{y^1}{F}\cfrac{\partial}{\partial 
x^1}+\cfrac{y^2}{F}\cfrac{\partial}{\partial 
x^2}=l^1\cfrac{\partial}{\partial 
x^1}+l^2\cfrac{\partial}{\partial x^2},
\end{split}
\end{equation*}
where $g$ is the determinant of the Hessian matrix $g_{ij}$.

The corresponding dual coframe is locally given by
\begin{equation*}
\begin{split}
\omega^1&=\cfrac{\sqrt g}{F}(y^2dx^1-y^1dx^2)=m_1dx^1+m_2dx^2\\
\omega^2&= \cfrac{\partial F}{\partial y^1}dx^1+\cfrac{\partial F}{\partial 
y^2}dx^2=l_1dx^1+l_2dx^2.
\end{split}
\end{equation*}

Next, one defines a moving coframing $(u;\omega^1,\omega^2,\omega^3)$ on $\pi^*TM$, orthonormal with respect to the
Riemannian metric on $\Sigma$ induced by the Finslerian metric $F$, where $u\in \Sigma$ and 
$\{\omega^1,\omega^2,\omega^3\}\in T^*\Sigma$. The
moving equations on this frame lead to the so-called Chern
connection. This is an almost metric compatible, torsion free connection of the
vector bundle $(\pi^*TM,\pi,\Sigma)$. 

Indeed, by a theorem of Cartan it follows that the coframe $(\omega^1,\omega^2,\omega^3)$ must satisfy the following structure equations
\begin{equation}\tag{2.2}
\begin{split}
d\omega^1&=-I\omega^1\wedge\omega^3+\omega^2\wedge\omega^3\\
d\omega^2&=-\omega^1\wedge\omega^3\\
d\omega^3&=K\omega^1\wedge\omega^2-J\omega^1\wedge\omega^3.
\end{split}
\end{equation}

The functions $I,J,K$ are smooth functions on $\Sigma$ called the {\bf invariants} of the Finsler structure $(M,F)$ in the sense of Cartan's equivalence problem (see for example \cite{BCS2000}, \cite{Br1997}, \cite{Br2002}).

This implies that on the vector bundle $\pi^*TM$ there exists a unique torsion-free and almost metric compatible connection $\nabla:C^\infty(T\Sigma)\otimes C^\infty(\pi^*TM)\to C^\infty(\pi^*TM)$, given by
\begin{equation}\tag{2.3}
 \nabla_{\hat{X}}Z=\{\hat{X}(z^i)+z^j\omega_j^{\ i}(\hat{X})\}e_i,
\end{equation}
where $\hat{X}$ is a vector field on $\Sigma$, $Z=z^ie_i$ is a section of $\pi^*TM$, and $\{e_i\}$ is the $g$-orthonormal frame field on $\pi^*TM$. 

The 1-forms $\omega_j^{\ i}$ define the {\bf Chern connection} of the Finsler structure $(M,F)$, where 
\begin{equation}\tag{2.4}\label{2.4}
(\omega_j^{\ i})=\left(\begin{array}{ccccc}{\omega_1{}^1} & {\omega_1{}^2}\\{\omega_2{}^1} 
&{\omega_2{}^2}\end{array}\right)=\left(\begin{array}{ccccc}{-I\omega^3} & 
{-\omega^3}\\{\omega^3} &{0}\end{array}\right),
\end{equation}
and $I:=A_{111}=A(e_1,e_1,e_1)$ is the {\bf Cartan scalar} for Finsler 
surfaces. Remark that $I=0$ is equivalent to the fact that the Finsler structure is 
Riemannian.

{\bf Remarks.}
\begin{enumerate}

\item We remark that the Chern connection gives a decomposition of the tangent bundle $T\Sigma$ by
\begin{equation*}
 T\Sigma=H\Sigma\oplus V\Sigma,
\end{equation*}
where the $H\Sigma$ is the horizontal distribution generated by $e_1,e_2$ and $V\Sigma$ is the vertical distribution generated by $\hat{e}_3$, where
$\hat{e}_1,\hat{e}_2,\hat{e}_3$ is the dual frame of the coframe $\omega^1,\omega^2,\omega^3$. 

\item For comparison, recall the structure equations of a Riemannian 
surface. They are obtained from (2.2) by putting $I=J=0$.

\item The scalar $K$ is called the {\bf Gauss curvature} of Finsler 
surface. In the case when $F$ is Riemannian, $K$ coincides with the usual 
 Gauss curvature of a Riemannian surface. 
\end{enumerate}

Differentiating again (2.2) one obtains the {\bf Bianchi identities}
\begin{equation}\tag{2.5}
\begin{split}
J&=I_2=\cfrac{1}{F}\biggl(y^1\cfrac{\delta I}{\delta x^1}+y^2\cfrac{\delta 
I}{\delta x^2}\biggr)\\
K_3&+KI+J_2=0,
\end{split}
\end{equation}
where 
$\{\frac{\delta}{\delta x^i},F\frac{\partial}{\partial y^i}\}$ is the adapted basis of $T\Sigma$, given by
\begin{equation*}
\frac{\delta}{\delta x^i}:=\frac{\partial}{\partial x^i}-N^j_i\frac{\partial}{\partial y^j}.
\end{equation*}
The functions $N_i^j$ are called the coefficients of the {\bf nonlinear connection} of $(M,F)$ (see \cite{BCS2000}, p. 33 for details).

The linear indices in $I_2$, $K_3$, $J_2$, etc. indicate differential terms with respect to 
$\omega_1,\omega_2,\omega_3$. For example 
$dK=K_1\omega^1+K_2\omega^2+K_3\omega^3$. The scalars $K_1$, $K_2$, $K_3$ 
are called the directional 
derivatives of $K$.

Nevertheless, remark that the scalars $I=I(x,y)$, $J=J(x,y)$, $K=K(x,y)$ and their derivatives 
live on $\Sigma$, not on $M$ as in the Riemannian case!

More generally, given any function $f:\Sigma \to \mathcal R$, one can write its 
differential in the form 
\begin{equation*}
df=f_1\omega_1+f_2\omega_2+f_3\omega_3.
\end{equation*}

Taking one more exterior differentiation of this formula, one obtains the 
following {\bf Ricci identities}:
\begin{equation}\tag{2.6}
\begin{split}
f_{21}-f_{12} &= -Kf_3 \\
f_{32}-f_{23} &= -f_1 \\
f_{31}-f_{13} &= If_1+f_2+Jf_3.
\end{split}
\end{equation}

One defines {\bf the curvature} of the Finsler structure $(M,F)$  as usual by
\begin{equation}\tag{2.7}
\Omega_j^{\ i}=d\omega_j^{\ i}-\omega_j^{\ k}\wedge \omega_k^{\ i},
\end{equation}
where $i,j,k\in\{1,2,3\}$, and $\omega_j^{\ i}$ is the Chern connection matrix (2.4). It easily follows that the only essential entry in the matrix 
$\Omega_j^{\ i}$ is 
\begin{equation}\tag{2.8}\label{2.8}
\Omega_2^{\ 1}=d\omega_2^{\ 1}=d\omega^{3}=K\omega^1\wedge\omega^2-J\omega^1\wedge\omega^3.
\end{equation}

We remark that the fact that the curvature 2-form is closed is a peculiarity of Finslerian surfaces that will be very useful in deriving the Gauss--Bonnet formula in the following sections.

Recall that a Finsler surface is called {\bf Landsberg} if the invariant $J$ vanishes. Bianchi identities imply that in this case $I_2=0$ and $K_3=-KI$. A Finsler structure having $I_1=0$, $I_2=0$ is called a {\bf Berwald} surface (see \cite{BCS2000}, Lemma 10.3.1, p. 267 for details). 

It is known that Berwald surfaces are in fact Riemannian surfaces if $K\neq 0$ or locally Minkowski flats if $K=0$ (see \cite{Sz1981} and
 \cite{BCS2000}, p. 278).

We also remark that on a Landsberg surface, even both $K$ and $g$
are quantities defined on the 3-dimensional manifold $\Sigma$, the product 
$K\sqrt{g}$ lives on $M$ (\cite{BCS2000}, p. 106).

\bigskip


Recall that the restriction of a Finsler norm to a tangent plane $T_xM$ gives a Minkowski norm on $T_xM$. For an arbitrary fixed $x\in M$, this Minkowski norm induces a Riemannian metric 
$\hat{g}$ on the punctured plane $\widetilde {T_xM}$ by
\begin{equation}\tag{2.9}
\hat{g}:=g_{ij}(y)dy^i\otimes dy^j,
\end{equation}
where $y=(y^i)$ are the global coordinates in $T_xM$.

 Remark that the Riemannian manifold $(\widetilde {T_xM}, \hat{g})$ is flat, i.e. the 
Gaussian curvature of $\hat{g}$ vanishes on $\widetilde {T_xM}$. This is a peculiarity of the two dimensional case (see \cite{BCS2000}, p. 388).

 The outward pointing normal to the indicatrix is 
\begin{equation}\tag{2.10}
\hat{n}_{out}=\cfrac{y}{F(y)}=\cfrac{y^i}{F(y)}\cdot\cfrac{\partial}{\partial 
y^i}.
\end{equation}

Indeed, let us consider $y^i=y^i(t)$ to be a unit speed parameterization of the indicatrix 
$\Sigma$. By derivation with respect to $t$ of the formula $g_{ij}(y)y^iy^j=1$ one obtains
$$
g_{ij}(y) y^i \dot y^j=0,
$$
where the dot notation means derivative with respect to $t$.

In the following let us consider the indicatrix $\Sigma_x$ as a Riemannian 
submanifold of the punctured Riemannian manifold $(\widetilde 
{T_xM},\hat{g})$, with the induced Riemannian metric $h$, and 
let $y(t)=(y^1(t),y^2(t))$ be a 
unit 
speed (with respect to $h$) parameterization of $\Sigma_x$.

Obviously, $F$ is Euclidean if and only if the main scalar $I$ restricted to $\Sigma_x$ vanishes. In other words, 
$I_{|_{\Sigma_x}}$ ``measures'' the deviation of $F$ on $T_xM$ from an Euclidean inner 
product.

The volume form of the Riemannian metric $\hat{g}$ on $T_xM$ is 
\begin{equation}\tag{2.11}
dV=\sqrt{g}dy^1 \wedge dy^2,
\end{equation}
where $ \sqrt{g}=\sqrt{det(g_{ij})}$, and the induced Riemannian volume 
form on the submanifold $\Sigma_x$ is
\begin{equation}\tag{2.12}
ds=\frac{\sqrt{g}}{F}(y^1\dot y^2-y^2\dot y^1)dt.
\end{equation}

Along $\Sigma_x$ the 1-form $ds$ coincides with
\begin{equation}\tag{2.13}
d\theta=\cfrac{\sqrt{g}}{F^2}(y^1dy^2-y^2dy^1).
\end{equation}
The parameter $\theta$ is called the {\bf Landsberg angle}.


{\bf Remarks.}

\begin{enumerate}
\item The formula $ds=\sqrt{g}(y^1\dot y^2-y^2\dot y^1)dt$ is valid as long 
as the underlying parameterization traces $\Sigma$ out in a positive manner.
\item The {\bf Riemannian length} of the indicatrix is therefore defined by
\begin{equation}\tag{2.14}
 L:=\int_{\Sigma_x} ds
\end{equation}
and it is typically NOT equal to $2\pi$ as in the case of Riemannian 
surfaces. This fact was remarked for the first time by M. Matsumoto 
\cite{M1986}. Since the indicatrix is a 1-dimensional submanifold its Riemannian length and the Riemannian volume are in fact identical. 
\end{enumerate}

\bigskip



The Riemannian length of the indicatrix $\Sigma_x$ is an integral where the 
integration domain also depends on $F$. One would like however to work with 
integrals over the {\it standard unit circle}
\begin{equation}\tag{2.15}
\mathbf S^1=\{y\in \widetilde{\mathbf T_xM}: (y^1)^2+(y^2)^2=1\},
\end{equation}
even with the price of a more complicated integrand.

It follows immediately that the indicatrix length in a Minkowski plane can be computed by 
\begin{equation}\tag{2.16}
L=\int_{\mathbf S^1}\cfrac{\sqrt{g}}{F^2}(y^1dy^2-y^2dy^1).
\end{equation}

\bigskip

Indeed, the 1-form
\begin{equation}\tag{2.17}
d\theta=\cfrac{\sqrt{g}}{F^2}(y^1dy^2-y^2dy^1)
\end{equation}
is a closed 1-form on $\widetilde{T_xM}$. By the use of Stokes' 
theorem one can easily see that integrating this over two corresponding arcs (see $\eqref{2.21}$ below) of $S$ and $\mathbf S^1$ 
one obtains the same answer (see \cite{BCS2000}, p. 101, 102).

One defines in this way the {\bf length function} of the indicatrix $\Sigma_x$ by
\begin{equation}\tag{2.18}
L:M\to [0,\infty),\qquad x\mapsto L(x)=\int_{\Sigma_x}\frac{1}{F}ds,
\end{equation}
or, equivalently,
\begin{equation}\tag{2.19}
  L(x)=\int_{\Sigma_x}d\theta.
\end{equation}

Let us also remark that 
\begin{equation}\tag{2.20}
 d\theta={\omega_2^{\ 1}}_{|_{\Sigma_x}},
\end{equation}
i.e. $d\theta$ is equal to the pure part $dy$ of $\omega_2^{\ 1}$,
therefore there is no harm if we write
\begin{equation}\tag{2.21}\label{2.21}
L(x)=\int_{\Sigma_x}\omega_2^{\ 1}.
\end{equation}
\bigskip

We define the {\bf Landsberg angle 
$\measuredangle _x(X,Y)$ of two Finslerian unit vectors} $X,Y\in \mathbf T_xM$ (with same 
origin, say $y=0$, or glided to have the same origin) and the tips on the indicatrix curve, as the oriented 
Riemannian angle of $X$ and $Y$ measured with the induced Riemannian metric 
$\hat{g}$.

In other words, for any two unit vectors $X,Y$ as above, their
Finslerian angle
is given by 
\begin{equation}\tag{2.22}
\measuredangle _x(X,Y):=\int_{\mathbf S_{\stackrel{\frown}{XY}}^1}
\cfrac{\sqrt{g}}{F^2}(y^1 d{y^2}-y^2 d{y^1})=
\int_{\Sigma_x\mid_{(X,Y)}}\sqrt{\hat{g}(\dot{y},\dot{y})}dt
=\int_{\Sigma_x\mid_{(X,Y)}}d\theta
,
\end{equation}
where ${\mathbf S_{\stackrel{\frown}{XY}}^1}$ and ${\Sigma_x\mid_{(X,Y)}}$ are the arcs on the unit Euclidean circle and the indicatrix curve described by the directions of the vectors $X$ and $Y$, respectively. Since the angle $\measuredangle _x(X,Y)$ is described by the integral of the 1-form $d\theta$, it is customary 
to call it the Landsberg angle.
\begin{center}
\includegraphics[height=7.5cm, angle=0]{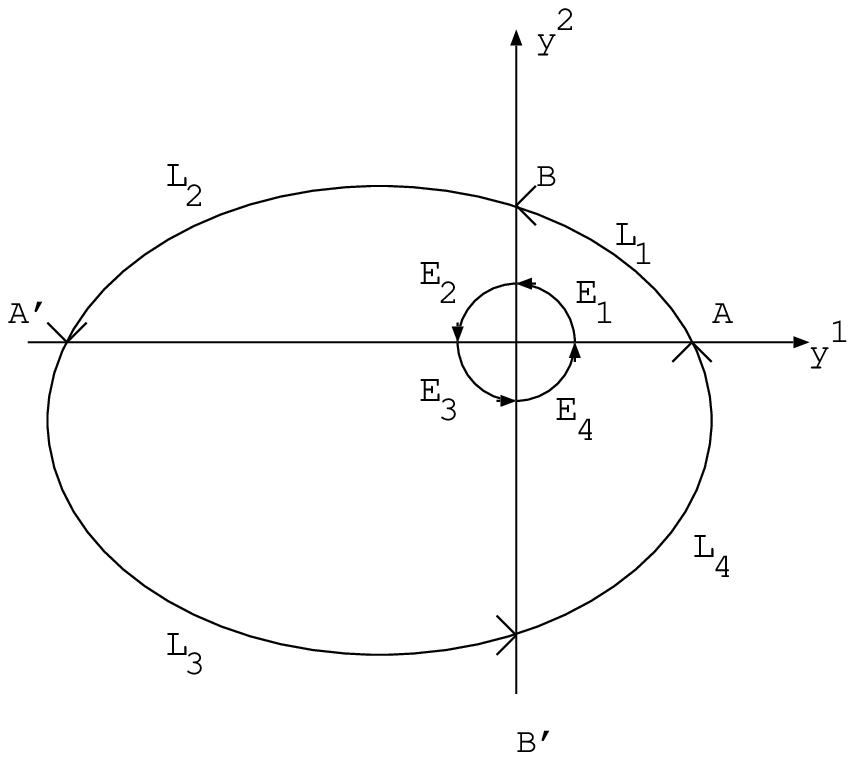}
\end{center}
{\bf Figure 1.} The Landsberg angle.

\bigskip

{\bf Remark.}

In this point it is important to remark that there are big differences between the Euclidean angles used in plane geometry and the Landsberg angles defined above (see Figure 1). Imagine the indicatrix of a Finsler space to be a translated ellipse (this is actually the case of a Randers metric) and the Euclidean unit circle inside it. We represented the Euclidean circle in the interior of the Finslerian indicatrix (they might actually intersect) only for making this explanation easy to follow. 
We denote the intersection points of the indicatrix with the coordinate axes by $A$, $B$, $A'$, $B'$, and the corresponding arcs by $L_1$, $L_2$, $L_3$, $L_4$, respectively. Moreover, we denote by  $E_1$, $E_2$, $E_3$, $E_4$ the corresponding arcs on the Euclidean unit circle. Obviously, the four Euclidean angles determined by the coordinate axes are all equal to $\cfrac{\pi}{2}$, and their sum equals $2\pi$. 

On the other hand, the Landsberg angles determined by the coordinate axes are described by the $\hat{g}$--Riemannian lengths of the indicatrix arcs $L_1$, $L_2$, $L_3$, $L_4$, respectively. One can easily see that the usual properties of angles known to hold good in Euclidean plane do not hold anymore. Indeed, remark for example that the opposite angle are not equal anymore,  $L_1\neq L_3$, $L_2\neq L_4$, nor the sum of adjacent angles equals $\pi$. However, we do know that the sum of $L_1$, $L_2$, $L_3$, $L_4$ equals the total length of indicatrix $L$.

A special case would be the case of an absolutely homogeneous Finsler norm, i.e. the case when the induced Minkowski norm satisfies the condition $F(-y)=F(y)$. In this case, the indicatrix, without being an ellipse, it is still a central symmetric curve, and therefore, the opposite angles are equal! In particular, 
$ L_1=L_3$ and $L_2=L_4$. 


\section{The normal lift of a curve}

$\quad$ Let us consider a smooth (or piecewise $C^\infty$) curve $\gamma:[0,r]\to M$ with the tangent vector $\dot{\gamma}(t)=T(t)$, parameterized such that $F(\gamma(t),\dot{\gamma}(t))=1$, and let $N$ be the {\bf  normal vector} along $\gamma$ defined by
 \begin{equation}\tag{3.1}\label{3.1}
\begin{split}
& g_N(N,N)=1 \\
& g_N(T,N)=0 \\
& g_N(T,T)=\sigma^2(t).
\end{split}
\end{equation}

We point out that here $g_N$ means
\begin{equation*}
(g_N)_{ij} = \frac{1}{2}\frac{\partial^2F^2}{\partial y^i\partial y^j}\big(\gamma(t), N(t)\big).
\end{equation*}

This kind of normal vector was
introduced by Z. Shen (\cite{Sh2001}, p. 27) and used by us in the formulation of Gauss--Bonnet
theorem for Landsberg surfaces with smooth boundary (\cite{SS2007}).

\bigskip

The {\bf normal lift} (shortly {\bf $N$-lift}) $\hat{\gamma}^\perp(t)$ of $\gamma(t)$ to $\Sigma$ is given by
\begin{equation}\tag{3.2}\label{3.2}
\begin{split}
& \hat{\gamma}^\perp:[0,r]\to \Sigma,\\
& \qquad t\mapsto \hat{\gamma}^\perp(t)=(\gamma(t),N(t)).
\end{split}
\end{equation}

\bigskip

The tangent vector $\hat{T}^\perp$ to the normal lift $\hat{\gamma}^\perp$ is given by
\begin{equation}\tag{3.3}\label{3.3}
\begin{split}
\hat{T}^\perp(t) & = \hat{\dot{\gamma}}^\perp(t)=\frac{d}{dt}\hat{\gamma}^\perp(t)=\dot{\gamma}^i(t)\frac{\partial}{\partial x^i}_{\vert_{(x,N)}}+\frac{d}{dt}N^i(t)\frac{\partial}{\partial y^i}_{\vert_{(x,N)}} \\
& =T^i(t)\frac{\delta}{\delta x^i}_{\vert_{(x,N)}}+(D_T^{(N)}N)^i\frac{\partial}{\partial y^i}_{\vert_{(x,N)}}.
\end{split}
\end{equation}

The local coefficients of the 
{\bf covariant derivative with reference vector $N$} along $\gamma$ are given by
\begin{equation}\tag{3.4}\label{3.4}
D_T^{(N)}U=(D_T^{(N)} U)^i\cdot\frac{\partial}{\partial x^i}_{\vert_{\gamma(t)}}=\Bigl[\frac{dU^i}{dt}+T^jU^k\Gamma^i_{jk}(x,N) \Bigr]\cdot\frac{\partial}{\partial x^i}_{\vert _{\gamma(t)}}
\end{equation}
for any $U=U^i(x)\frac{\partial}{\partial x^i}$ vector field along $\gamma$, where $\Gamma^i_{jk}$ are the Chern connection coefficients, i.e. $\omega_i^j=\Gamma_{ik}^jdx^k$.

Remark that the term $D^{(N)}_TN$ in $\eqref{3.3}$ means the covariant derivative of $N$ along $\gamma$ with reference vector $N$.\\

We recall here a useful lemma (Lemma 7.1 in \cite{SS2007}):

{\bf Lemma 3.1.}
\begin{equation}\tag{3.5}\label{3.5}
\frac{d}{dt}g_N(V,W)=g_N(D_T^{(N)}V,W)+g_N(V,D_T^{(N)}W)+
2A(V,W,D_T^{(N)}N)_{\vert_{(x,N)}},
\end{equation}
where $A$ is the Cartan tensor (see \cite{BCS2000}, p. 30). 

Using this, we obtain 
\begin{equation}\tag{3.6}\label{3.6}
\begin{split}
& g_N(D_T^{(N)}N,N)=0,\\
& g_N(D_T^{(N)}T,N)+g_N(T, D_T^{(N)}N)=0,\\
& g_N(D_T^{(N)}T,T)=
\sigma(t)\frac{d\sigma}{dt}-A(T,T,D_T^{(N)}N).
\end{split}
\end{equation}

Similarly to the notion of Finslerian geodesics we can define the notion of $N$-parallel of a Finsler structure.\\
\newline
{\bf Definition 3.1.}$\;$
A curve $\gamma$ on the surface $M$, in Finslerian natural parameterization,
 is called an
{\bf N-parallel} of the Finslerian structure $(M,F)$ if and only if we have
\begin{equation}\tag{3.7}\label{3.7}
D_T^{(N)}N=0.
\end{equation}

It follows from $\eqref{3.3}$ that the tangent vector to the normal lift of an  $N$-parallel curve $\gamma$
on $M$ is given by 
$\hat{T}^\perp=T^i\frac{\delta}{\delta x^i}_{\vert_{(x,N)}}$. In other words, we obtain the following characterization of an $N$-parallel of a Finsler surface.\\
\newline 
{\bf Proposition 3.1.}$\;$ {\it A curve on $M$ is an  $N$-parallel curve if and only if its
normal vector $N$ is transported parallel along $\gamma$.} 

\bigskip

{\bf Remarks.}
\begin{enumerate}
\item If $\gamma$ is an  $N$-parallel, then we have
\begin{equation}\tag{3.8}\label{3.8}
\begin{split}
& g_N(D_T^{(N)}T,N)=0\\
& g_N(D_T^{(N)}T,T)=
\sigma(t)\frac{d\sigma}{dt}.
\end{split}
\end{equation}
\item The curve $\gamma$ is an  $N$-parallel if and only if 
$\nabla_{\hat{T}^\perp}l=D_T^{(N)}N=0$. This implies that
      $g_N(\nabla_{\hat{T}^\perp}l,l)=0$, i.e. $\nabla_{\hat{T}^\perp}l$
      is orthogonal to the indicatrix.
\end{enumerate}

\bigskip

In case of an arbitrary curve $\gamma$ on $M$, from 
\begin{equation*}
g_N(D_T^{(N)}N,N)=0, \qquad g_N(T,N)=0,
\end{equation*}
it follows that the vector $D_T^{(N)}N$ is proportional to $T$,
i.e. there exists a non-vanishing function $k_{T}^{(N)}(t)$ such that
\begin{equation}\tag{3.9}\label{3.9}
D_T^{(N)}N=-\frac{k_{T}^{(N)}(t)}{\sigma^2(t)}T,\qquad \sigma(t)\neq 0.
\end{equation}

The function $k_{T}^{(N)}(t)$ will be called the {\bf $N$-parallel curvature
of $\gamma$}. The minus sign is put only in order to obtain the
same formulas as in the classical theory of Riemannian manifolds.

In other words, we have
\begin{equation}\tag{3.10}\label{3.10}
g_N(D_T^{(N)}N,T)=-g_N(D_T^{(N)}T,N)=-k_{T}^{(N)}(t).
\end{equation}

Since $\{N,T\}$ is a basis, we also obtain
\begin{equation}\tag{3.11}\label{3.11}
D_T^{(N)}T=k_{T}^{(N)}(t)N+B(t)T,
\end{equation}
where we put 
 \begin{equation}\tag{3.12}\label{3.12}
B(t)=\frac{1}{\sigma(t)}\frac{d\sigma(t)}{dt}-\frac{1}{\sigma^2(t)}A_{\vert_{(x,N)}}(T,T,D_T^{(N)}N).
\end{equation}

By making use of the cotangent map of $\hat{\gamma}^{\perp}$ we compute
\begin{equation}\tag{3.13}\label{3.13}
\begin{split}
& \hat{\gamma}^{\perp *}\omega^1\frac{\partial}{\partial
 t}=\omega^1(\hat{T}^\perp)_{\vert_{(x,N)}}
=\frac{\sqrt{g}}{F}(N^2T^1-T^2N^1)=\sigma(t)\\
& \hat{\gamma}^{\perp *}\omega^2\frac{\partial}{\partial
 t}=\omega^2(\hat{T}^\perp)_{\vert_{(x,N)}}=g_N(N,T)=0\\
& \hat{\gamma}^{\perp *}\omega^3\frac{\partial}{\partial
 t}=\omega^3(\hat{T}^\perp)_{\vert_{(x,N)}}
=\frac{\sqrt{g}}{F}
\Bigl[ N^2(D_T^{(N)}N)^1-N^1(D_T^{(N)}N)^2
\Bigr]=-\frac{k_{T}^{(N)}(t)}{\sigma(t)}
\end{split}
\end{equation}
(for details see also \cite{SS2007}).

Therefore we obtain
\begin{equation}\tag{3.14}\label{3.14}
\begin{split}
& \hat{\gamma}^{\perp *}\omega ^1=\sigma(t)dt \\
& \hat{\gamma}^{\perp *}\omega^2=0 \\
& \hat{\gamma}^{\perp *}\omega^3=-\frac{k^{(N)}_{T}}{\sigma(t)}dt.
\end{split}
\end{equation}

If we denote by $\{\hat e_1, \hat e_2, \hat e_3\}$ the dual frame on $\Sigma$ of the orthonormal coframe $\{\omega^1, \omega^2, \omega^3\}$, then we obtain that the tangent vector to the normal
lift of $\hat{\gamma}^{\perp}$ is 
\begin{equation}\tag{3.15}\label{3.15}
\hat{T}^\perp=\sigma(t)\hat{e}_1-\frac{k^{(N)}_{T}}{\sigma(t)}\ \hat{e}_3\in<\hat{e}_1,\hat{e}_3>,
\end{equation}
where $<\hat{e}_1,\hat{e}_3>$ is the 2-plane generated by $\hat{e}_1$, $\hat{e}_3$.
\bigskip


Remark that in the case when $\gamma$ is an $N$-parallel, we have
\begin{equation}\tag{3.16}\label{3.16}
D_T^{(N)}T=\frac{1}{\sigma(t)}\frac{d\sigma(t)}{dt}T,
\end{equation}
and
\begin{equation}\tag{3.17}\label{3.17}
\begin{split}
& \hat{\gamma}^{\perp *} \omega ^1=\sigma(t)dt \\
& \hat{\gamma}^{\perp *}\omega^2=0 \\
& \hat{\gamma}^{\perp *}\omega^3=0.
\end{split}
\end{equation}

 Finally, we remark that the tangent vector to the normal
lift of an $N$-parallel is 
\begin{equation}\tag{3.18}\label{3.18}
\hat{T}^\perp=\sigma(t)\hat{e}_1.
\end{equation}

\bigskip


{\bf Remark.}

Let us remark that the $N$-lift used in this section is different from the {\bf canonical lift} (or {\bf tangential lift}) of a curve. 
Indeed, for an arbitrary curve $\gamma:[0,r]\to M$ with the usual properties, 
the {\bf canonical lift}  of $\gamma$ to $\Sigma$ is given by
\begin{equation}\tag{3.19}\label{3.19}
\begin{split}
& \hat{\gamma}:[0,r]\to \Sigma,\\
& \qquad t\mapsto \hat{\gamma}(t)=(\gamma(t),\dot{\gamma}(t)).
\end{split}
\end{equation}
This is well defined because $(\gamma(t),\dot{\gamma}(t))\in \Sigma$ because of the Finslerian natural parameterization of $\gamma$ 
(see \cite{BCS2000}, p. 112).

Let us consider now the normal vector $\mathcal V$ along $\gamma$ with respect to the tangent vector $T$ defined by
\begin{equation*}
g_T(T,\mathcal{V})=0.
\end{equation*}

Here, by $g_T$ we mean
\begin{equation*}
\big(g_T\big)_{ij} = \frac{1}{2}\frac{\partial^2F^2}{\partial y^i\partial y^j}\big(\gamma(t), T(t)\big).
\end{equation*}

We have the fundamental relations
\begin{equation*}
\begin{split}
& g_T(T,T)=1 \\
& g_T(T,\mathcal{V})=0,
\end{split}
\end{equation*}
and let us put
\begin{equation*}
 \mu^2(t):=g_T(\mathcal{V},\mathcal{V}).
\end{equation*}

We also obtain 
\begin{equation}\tag{3.20}\label{3.20}
\begin{split}
& g_T(D_T^{(T)}T,T)=0\\
& g_T(D_T^{(T)}T,\mathcal{V})+g_T(D_T^{(T)}\mathcal{V},T)=0\\
& g_T(D_T^{(T)}\mathcal{V},\mathcal{V})=
\mu(t)\frac{d\mu}{dt}-A(\mathcal{V},\mathcal{V},D_T^{(T)}T).
\end{split}
\end{equation}
The local coefficients of the 
{\bf covariant derivative with reference vector $T$} along $\gamma$ are given by
\begin{equation}\tag{3.21}\label{3.21}
D_T^{(T)}U=(D_T^{(T)} U)^i\cdot\frac{\partial}{\partial x^i}_{\vert_{\gamma(t)}}=\Bigl[\frac{dU^i}{dt}+T^jU^k\Gamma^i_{jk}(x,T) \Bigr]\cdot\frac{\partial}{\partial x^i}_{\vert _{\gamma(t)}}
\end{equation}
for any $U = U^i(x)\frac{\partial}{\partial x^i}$ vector field along $\gamma$, where $\Gamma^i_{jk}$ are the Chern connection local coefficients.

One can see that the term $D_T^{(T)}T$ in $\eqref{3.20}$ means the covariant derivative of $T$ along $\gamma$ with reference vector $T$.

From 
\begin{equation*}
g_T(D_T^{(T)}T,T)=0, \qquad g_T(T,\mathcal{V})=0,
\end{equation*}
it follows that the vector $D_T^{(T)}T$ is proportional to
$\mathcal{V}$, i.e. there exists a non-vanishing function $k_{\mathcal{V}}^{(T)}(t)$ such that
\begin{equation*}
D_T^{(T)}T=k_{\mathcal{V}}^{(T)}(t)\mathcal{V}.
\end{equation*}
The function $k_{\mathcal{V}}^{(T)}(t)$ is called the {\bf signed curvature of $\gamma$ over T}.

On the other hand, remark that from 
$g_T(D_T^{(T)}T,\mathcal{V})+g_T(D_T^{(T)}\mathcal{V},T)=0$
we obtain 
$g_T(D_T^{(T)}\mathcal{V},T)=-g_T(D_T^{(T)}T,\mathcal{V})=-k_{\mathcal{V}}^{(T)}(t)$.

We also obtain
\begin{equation}\tag{3.22}\label{3.22}
\begin{split}
& \hat{\gamma}^* \omega ^1=0 \\
& \hat{\gamma}^*\omega^2=dt \\
& \hat{\gamma}^*\omega^3=-\frac{k^{(T)}_{\mathcal{V}}}{\mu}dt,\qquad (\mu\neq 0).
\end{split}
\end{equation}

In the case when $\gamma$ is a Finslerian geodesic, we have by definition $D_T^{(T)}T=0$
and therefore 
\begin{equation*}
\begin{split}
& g_T(T,D_T^{(T)}\mathcal{V})=0\\
& g_T(\mathcal{V},D_T^{(T)}\mathcal{V})=\mu(t)\frac{d\mu(t)}{dt}.
\end{split}
\end{equation*}
It follows
\begin{equation*}
D_T^{(T)}\mathcal{V}=\frac{1}{\mu(t)}\frac{d\mu(t)}{dt}\mathcal{V}
\end{equation*}
and
\begin{equation}\tag{3.23}\label{3.23}
\begin{split}
& \hat{\gamma}^* \omega ^1=0 \\
& \hat{\gamma}^*\omega^2=dt \\
& \hat{\gamma}^*\omega^3=0.
\end{split}
\end{equation}

 The tangent vector to the tangential
lift $\hat\gamma$ of a Finslerian geodesic $\gamma$ is 
\begin{equation*}
\hat{T}=\hat{e}_2.
\end{equation*}

\bigskip


We will end this section by pointing out that this theory reduces to the
classical theory in the case of a {\bf  Riemannian surface}. 

Let us assume that our Finslerian structure on $M$ is actually a
Riemannian one, and let us denote the Riemannian metric on the surface
$M$ by $g$. Then the normal along a curve $\gamma$ on $M$, natural
parameterized,  is defined
by 
\begin{equation}\tag{3.24}\label{3.24}
\begin{split}
& g(T,T)=1\\
& g(T,N)=0\\
& g(N,N)=1.
\end{split}
\end{equation}
Therefore the two types of normals $N$ and $\mathcal{V}$ defined above
coincide, and $\sigma=\mu=1$.

The tangent lift of $\gamma$ to $\Sigma^g$ is 
\begin{equation*}
\begin{split}
& \hat{\gamma}:[0,r]\to \Sigma^g,\\
& \qquad t\mapsto \hat{\gamma}(t)=(\gamma(t),T(t)),
\end{split}
\end{equation*}
where $\Sigma^g$ is the total space of the unit sphere bundle of the Riemannian structure $(M,F)$. Its tangent vector is
\begin{equation*}
\begin{split}
\hat{T}(t) & = \hat{\dot{\gamma}}(t)=\frac{d}{dt}\hat{\gamma}(t)=\dot{\gamma}^i(t)\frac{\partial}{\partial x^i}_{\vert_{(x,T)}}+\frac{d}{dt}T^i(t)\frac{\partial}{\partial y^i}_{\vert_{(x,T)}} \\
& =T^i(t)\frac{\delta}{\delta x^i}_{\vert_{(x,T)}}+(D_TT)^i\frac{\partial}{\partial y^i}_{\vert_{(x,T)}},
\end{split}
\end{equation*}
where $D_TT$
is the usual covariant derivative along $\gamma$ with respect to the Levi Civita
connection of $g$.

By derivation we obtain
\begin{equation*}
\begin{split}
& g(D_TT,T)=0\\
& g(D_TT,N)+g(D_TN,T)=0\\
& g(D_TN,N)=0.
\end{split}
\end{equation*}

From $g(D_TT,T)=0$, $g(T,N)=0$ it follows
\begin{equation*}
D_TT=k_N(t)N,
\end{equation*}
where the function $k_N(t)$ is the usual Riemannian signed curvature of
$\gamma$. 

On the other hand, from 
\begin{equation*}
\begin{split}
& g(T,D_TN)=-g(D_TT,N)=-k_N(t)\\
& g(D_TN,N)=0
\end{split}
\end{equation*}
we obtain
\begin{equation*}
D_TN=-k_N(t)T.
\end{equation*}

Let us consider now the $N$-lift of $\gamma$ to $\Sigma^g$ defined as
above. By similar computations as in the Finslerian case, in the
Riemannian case we obtain 
\begin{equation*}
k_T^{(N)}=k_N=k_N^{(T)},
\end{equation*}
i.e. the $N$-parallel curvature and the signed curvature over $T$ coincide with the usual
Riemannian signed curvature.

Moreover, the curve $\gamma$ is a Riemannian geodesic if and only if one
of the following relations hold
\begin{enumerate}
\item $k_N=0$
\item $D_TT=0$
\item $D_TN=0$,
\end{enumerate}
i.e. on a Riemannian geodesic the vectors $T$ and $N$ are equally
parallel transported along $\gamma$. In other words, on a Riemannian
manifold, the Riemannian geodesics and the $N$-parallel curves coincide.




\section{The Gauss--Bonnet theorem for Finsler surfaces with smooth boundary}

$\quad$ The proof of the Gauss-Bonnet theorem for Finsler manifolds without boundary was given by D. Bao and S. S. Chern in \cite{BC1996} using the transgression method. Using their method we have extended the result to Landsberg surfaces with smooth boundary \cite{SS2007}. 

In the present paper, 
we are going to give a Gauss-Bonnet type formula for Riemann--Finsler surfaces where the indicatrix volume does not need to be constant anymore, using an idea of B. Lackey \cite{L2002}. 

\bigskip

We start by discussing the case of a Riemann--Finsler surface with {\bf smooth boundary}.

Let $(M,F)$ be a compact Riemann--Finsler surface 
and $D\subset M$ a domain with smooth  boundary 
 $\partial D=\gamma:[a,b]\mapsto M$, given by $x^i=x^i(t)$. We assume $\gamma$ to be unit speed, i.e. $F(\gamma(t),T(t))=1$, where $T(t)=\dot{\gamma}(t)$.

\bigskip

{\bf Proposition 4.1.}

{\it 
Let $(M,F)$ be a compact oriented Finslerian surface and $D\subset M$ a 
domain with boundary $\partial D$.  Let $N:\partial D\to \Sigma$ be the inward pointing Finslerian unit normal on $\partial D$. 

Then, we have
\begin{equation}\tag{4.1}\label{4.1}
\begin{split}
 \int_D \cfrac{1}{L(x)}\Bigl[X^*(K\ \omega^1\wedge \omega^2
 & - J\omega^1\wedge\omega^3)-d\log L(x)\wedge X^*(\omega^3)
\Bigr]\\
& +\int_{N(\partial D)}\cfrac{1}{L(x)}\ \omega_1^{\ 2}=
\mathcal X(D),
\end{split}
\end{equation}
where 
$L(x)$ is the Riemannian length of the indicatrix $\Sigma_x$,
$X$ is a unit prolongation of $N$, 
$K$ is the Gauss curvature, 
and $\mathcal X(D)$ is the Euler characteristic of $D$. 
}

\bigskip
 
The proof follows \cite{BCS2000}, p. 106 or \cite{SS2007}.
Indeed, remark first that we can extend the normal vector field $N$ on $\gamma$ to a vector field $V$ on $M$ with only finitely many zeros 
$x_1, x_2, \dots, x_k$ in $D\setminus \partial D$. It is then known that the sum of indices of $X$ is equal to the Euler characteristic $\mathcal X(D)$ (see for example \cite{Spiv1979}, Vol. V, p. 561).

By removing from $D$ the interiors of the geodesic circles 
$S_{\alpha}^\varepsilon$ 
(centered at $x_{\alpha}$ of radius $\varepsilon >0$), one obtains the 
manifold with boundary $D_{\varepsilon}$. 
Remark that in this case, the boundary of $D_{\varepsilon}$ consists of 
the boundaries of the geodesic circles $S_{\alpha}^\varepsilon$ and the boundary of $D$. 

Since $V$ has all zeros in $D\setminus \partial D$, it follows that $V$ has no zeros on 
$D_{\varepsilon}$ and therefore we can normalize it obtaining in this way the application
\begin{equation}\tag{4.2}\label{4.2}
X=\cfrac{V}{F(V)}:D_{\varepsilon}\to \Sigma,\quad  x\mapsto
 \cfrac{V(x)}{F(V(x))}\ .
\end{equation}

Using $X$ we can lift $D_{\varepsilon}$ to $\Sigma$ constructing in this way the 2-dimensional submanifold $X(D_{\varepsilon})$ of $\Sigma$ such that we can integrate formula $\eqref{2.8}$ over this submanifold. 

\bigskip

However, before doing this, we make the following remark. 

From the degree theory (see for example \cite{Mil1965}) it results that
\begin{equation*}
\lim_{\varepsilon\to 0}\int_{X(S^\varepsilon _\alpha)}\omega_1{}^2= -i_{\alpha}(X)
\int_{\Sigma_{x_{\alpha}}}\omega_1{}^2=-i_{\alpha}(X)L(x_\alpha),
\end{equation*}
where $i_{\alpha}(X)$ is the index of $X$ at $x_{\alpha}$. Here the indicatrix $\Sigma_{x_{\alpha}}$ is traced in the 
counterclockwise orientation.

Since all the indicatrices are smooth closed convex curves inclosing the origin, it follows that $L(x)\neq 0$ for any $x\in M$. Using Lackey's idea (\cite{L2002}),we compute the index of $X$ at an arbitrary fixed zero point 
$x_\alpha$ by
\begin{equation}\tag{4.3}\label{4.3}
 -i_{\alpha}(X)=\cfrac{1}{L(x_\alpha)}\ \lim_{\varepsilon \to 0}\int_{X(S^\varepsilon _\alpha)}\omega_1{}^2=
\ \lim_{\varepsilon \to 0}\int_{X(S^\varepsilon _\alpha)}\cfrac{1}{L(x_\alpha)}\ \omega_1{}^2,
\end{equation}
where we have used the fact that when taking the limit of the integral 
$\omega_1{}^2$ the $x$ terms actually  do not contribute anymore 
because 
the metric radius continuously shrinks.

By summing over the zeros of $X$ and using Stokes' Theorem, it follows
\begin{equation}\tag{4.4}\label{4.4}
 \begin{split}
  -\mathcal X(D) & = -\sum_{\alpha=1}^k i_\alpha (X)=\sum_{\alpha=1}^k 
\lim_{\varepsilon \to 0}\int_{X(S^\varepsilon _\alpha)}\cfrac{1}{L(x_\alpha)}\ \omega_1{}^2
= \sum_{\alpha=1}^k 
\lim_{\varepsilon \to 0}\int_{X(S^\varepsilon _\alpha)}\Pi\\
& =\sum_{\alpha=1}^k 
\lim_{\varepsilon \to 0}\int_{X(S^\varepsilon _\alpha)}\Pi
+\int_{N(\partial D)}\Pi-\int_{N(\partial D)}\Pi=\int_{\partial X(D_\varepsilon)}\Pi-\int_{N(\partial D)}\Pi\\
& = \int_DX^*(d\Pi)-\int_{N(\partial D)}\Pi,
 \end{split}
\end{equation}
where we put
\begin{equation}\tag{4.5}\label{4.5}
 \Pi:=\frac{1}{L(x)}\ \omega_1{}^2.
\end{equation}

In this way we obtain the following

\bigskip

{\bf Topological Lemma.}

{\it 
Let $(M,F)$ be a compact oriented Finslerian surface and $D\subset M$ a domain with smooth boundary $\partial D$. Let $N:\partial D\to \Sigma$ be the inward pointing Finslerian unit normal on $\partial D$. 

Then, we have
\begin{equation}\tag{4.6}\label{4.6}
 -\int_DX^*(d\Pi)+\int_{N(\partial D)}\Pi= \mathcal X(D),
\end{equation}
where the notations are the same as above.
}

This is the extension of the topological lemma in \cite{L2002} to the case of Finsler surfaces with smooth boundary.
\bigskip

We need now to compute the first term of the sum in left hand side of (\ref{4.6}). 

In order to prove $\eqref{4.1}$, a straightforward computation gives
\begin{equation*}
\begin{split}
 d\Pi & =d\Bigl[ \frac{1}{L(x)}\ \omega_1{}^2\Bigr]=
\frac{1}{L(x)}\Bigl[ d\omega_1{}^2-d\log L(x)\wedge \omega_1{}^2  \Bigr]\\
 & = \frac{1}{L(x)}\Bigl[ -K\omega^1\wedge \omega^2+J\omega^1\wedge \omega^3-d\log L(x)\wedge \omega_1{}^2  \Bigr],
\end{split}
\end{equation*}
and Proposition 4.1 is proved.

Let us remark that in the first integral of the sum in the left hand side
 of (4.1) we should have written $L\circ X(x)$ instead of $L(x)$. However, since 
$L(x)$ is the length of the indicatrix at $x\in M$, and $X$ is an unit vector field on $M$, one can easily see that $L(x)$ and   $L\circ X(x)$ give actualy the 
same value.
The same is true for the unit vector $N$ and we will simplify the notation in this paper by writing simply $L(x)$.

\bigskip

We can evaluate the second term in the left hand sum of $\eqref{4.1}$ as follows:
\begin{equation}\tag{4.7}\label{4.7}
\begin{split}
 \int_{N(\partial D)}\frac{1}{L(x)}\ \omega_1^2 & = 
\int_{\partial D}N^*\Bigl(\frac{1}{L(x)}\ \omega_1^2\Bigr)=
\int_{\partial D}\frac{1}{L(x)}\ N^*\Bigl(\omega_1^2\Bigr)\\
& =\int_{\gamma}\frac{1}{L(x)}\frac{k_T^{(N)}(t)}{\sigma(t)}dt,
\end{split}
\end{equation}
where we have used (\ref{3.14}).

From Proposition 4.1 and formula (\ref{4.7}) we conclude

\bigskip
{\bf Theorem 4.2.} (The Gauss--Bonnet formula for Finsler surfaces with smooth boundary)

{\it Let $(M,F)$ be a compact oriented Finslerian surface and $D\subset M$ a domain with boundary $\partial D=\gamma$. Let $N:\partial D\to \Sigma$ be the inward pointing Finslerian unit normal on $\partial D$. 

Then, we have
\begin{equation}\tag{4.8}\label{4.8}
\begin{split}
 \int_D \cfrac{1}{L(x)}\Bigl[X^*(K\ \omega^1\wedge \omega^2
 & - J\omega^1\wedge\omega^3)-d\log L(x)\wedge X^*(\omega^3)
\Bigr]\\
& +\int_{\gamma}\frac{1}{L(x)}\frac{k_T^{(N)}(t)}{\sigma(t)}dt=
\mathcal X(D),
\end{split}
\end{equation}
where 
$L(x)$ is the Riemannian length of the indicatrix $\Sigma_x$,
$N$ is the inward pointing normal to the boundary $\partial D$,
$X$ is a unit prolongation of $N$, 
$K$ is the Gauss curvature, 
and $\mathcal X(D)$ is the Euler characteristic of $D$. 
}


{\bf Remarks.}
\begin{enumerate}
 \item If $(M,F)$ is a Landsberg surface, 
then $J=0$, $L(x)=L=\text{constant}$ and therefore (\ref{4.8}) gives the Gauss--Bonnet formula for Landsberg surfaces (see \cite{BC1996}, \cite{BCS2000} for the boundaryless case and \cite{SS2007} for the smooth boundary case). In other words we have
\begin{equation}\tag{4.9}\label{4.9}
 \cfrac{1}{L}\int_D K\sqrt{g} dx^1\wedge dx^2
 +\cfrac{1}{L}\int_{\gamma}\frac{k_T^{(N)}(t)}{\sigma(t)}dt=
\mathcal X(D),
\end{equation}
with the same notations as above and where $g$ is the determinant of the induced metric $g_{ij}$.

\item If $M$ is a compact orientable boundaryless Finsler manifold, then the Gauss-Bonnet formula reads
\begin{equation}\tag{4.10}\label{4.10}
 \int_D \cfrac{1}{L(x)}\Bigl[X^*(K\ \omega^1\wedge \omega^2
  - J\omega^1\wedge\omega^3)-d\log L(x)\wedge X^*(\omega^3)
\Bigr]=
\mathcal X(D).
\end{equation}
One can see that this formula agrees with \cite{Sh1996}.

\item If $(M,F)$ is Riemannian, then one obtains immediately the usual Gauss--Bonnet formula for Riemannian surfaces with smooth boundary (see for example \cite{Spiv1979}, p. 558, \cite{SST2003}, p. 34 and many other places).

\end{enumerate}

\bigskip


\section{The Gauss--Bonnet theorem for Finsler surfaces with regular piecewise $C^\infty$ boundary}

$\quad$ Let $(M,F)$ a compact Finsler  surface and $D\in M$ a domain with regular piecewise
$C^\infty$ boundary $\partial D=\gamma:[a,b]\mapsto M$, given by
$x^i=x^i(t)$. Let $a=t_0<t_1<\dots <t_k=b$ be a partition of $[a,b]$
such that $\gamma$ is $C^\infty$ on each closed subinterval
$[t_{s-1},t_s]$, $s\in\{1,2,\dots, k\}$. We assume $\gamma$ to be unit speed, i.e. $F(\gamma(t),T(t))=1$, where $T(t)=\dot{\gamma}(t)$.

For the sake of simplicity, let us assume that our boundary curve 
$\gamma$ has only one corner $x_0=x(t_0)$, for some $t_0 \in [a,b]$. In the case of $k$ corners, we are going to sum the quantities to be obtained below.

As in the proof of Theorem 4.2., we take the $N$-lift of $\gamma$ to $\Sigma$:
\begin{equation}\tag{5.1}\label{5.1}
\hat{\gamma}^\perp:[a,b]\to \Sigma,\qquad t\mapsto (x(t),N(t)),
\end{equation}
where $N$ is defined as above by $g_{N(t)}(T(t),N(t))=0$ for all 
$t\in [a,b]\setminus\{t_0\}$. 

Remark that in the case of one corner, the $N$-lift $\hat{\gamma}^\perp$
is not a closed curve anymore.

Indeed, let us denote by $T^-$ and $T^+$ the tangent vectors to $\gamma$ in $x_0$, i.e.
\begin{equation}\tag{5.2}\label{5.2}
T^-=\lim_{t\nearrow t_0}T(t),\qquad T^+=\lim_{t\searrow t_0}T(t),
\end{equation}
and define the corresponding normals at $x_0$ by
\begin{equation}\tag{5.3}\label{5.3}
g_{N^-}(N^-,T^-)=0,\qquad g_{N^+}(N^+,T^+)=0,
\end{equation}
respectively.

\smallskip
\begin{center}
\includegraphics[height=8cm, angle=0]{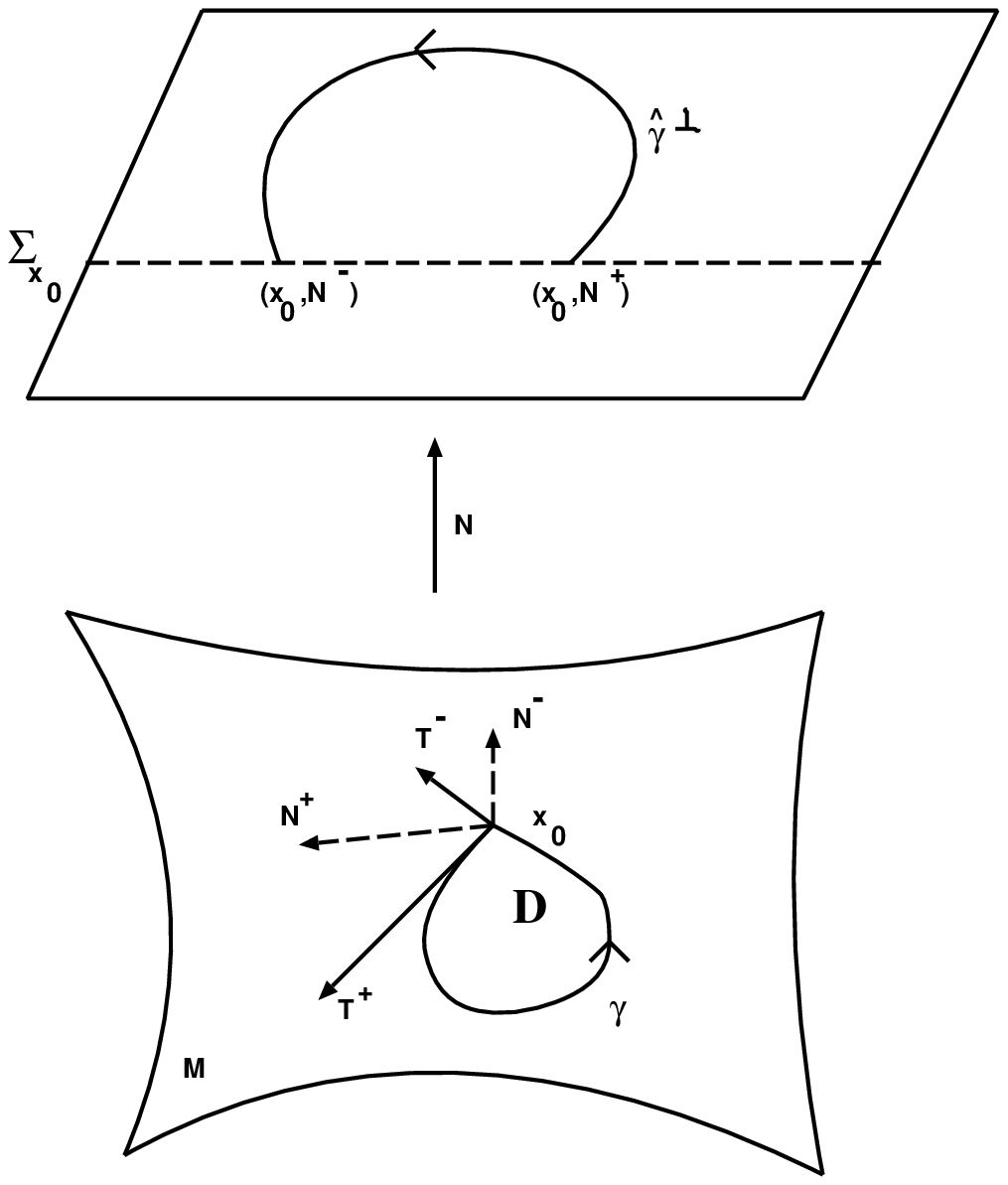}
\end{center}
{\bf Figure 2.} The normal lift of a regular $C^\infty$ piecewise curve with a corner.
\smallskip

It follows that at the point $x_0$ the tangent vector $T(t)$ has a discontinuous jump from $T^-$ to $T^+$, and similarly, the normal vector $N(t)$ has also a discontinuous jump from $N^-$ to $N^+$.

When lifting the curve $\gamma$ to $\Sigma$ we obtain a $C^{\infty}$ curve $\hat{\gamma}$ in $\Sigma$ with the ends $(x_0,N^-)$ and $(x_0,N^+)$. Remark that 
$N^-,N^+\in \widetilde{T_{x_0}M}$, and $F(x_0,N^-)=F(x_0,N^+)=1$ (see Figure 2).

Now, since $N^-$ and $N^+$ are two vectors in $T_{x_0}M$ with the origin in $x_0$ and the tips on the indicatrix, their Landsberg angle is
\begin{equation}\tag{5.4}\label{5.4}
\begin{split}
\measuredangle _{x_0} (N^-,N^+) & =\int_{{\stackrel{\frown}{N^-N^+}}}d\theta
=\int_{{\stackrel{\frown}{N^-N^+}}}
\cfrac{\sqrt{g}}{F^2}(y^1 d{y^2}-y^2 d{y^1})\\
 & =\int_{\tau_1}^{\tau_2}\sqrt{\hat{g}(\dot{y},\dot{y})}d\tau,
\end{split}
\end{equation}
where $y=y(\tau)$ is a unit speed parameterization of the indicatrix and 
$N^-=y(\tau_1)$, $N^+=y(\tau_2)$. Here the Landsberg angle is always evaluated using the positive orinted indicatrix arc joining the points $N^-$, $N^+$. Here the positive orientation on the indicatrix is given by $ds$.

We will proceed further and extend the normal vector field $N$ along $\gamma$ to a smooth section of $TM$ defined along the subset $\gamma\subset M$.

\smallskip
\begin{center}
\includegraphics[height=5cm,width=8cm, angle=0]
{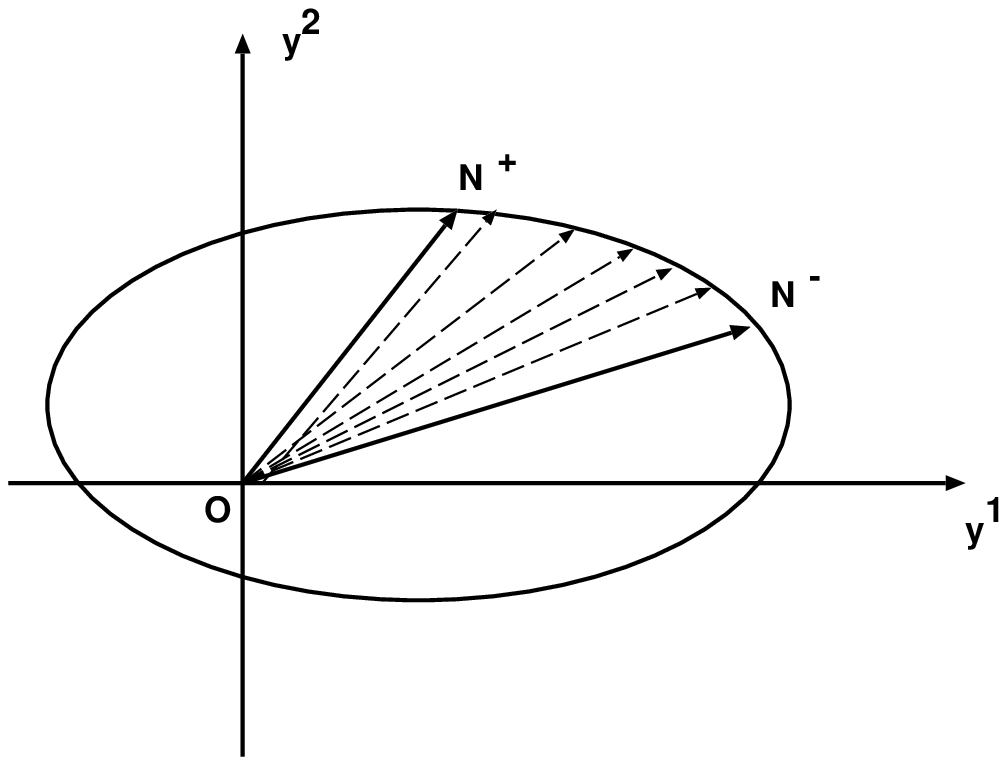}
\end{center}
{\bf Figure 3.} The Landsberg angle on a Finslerian indicatrix.

\smallskip

Intuitively, the most natural way of doing this is to consider the set of vectors in $T_{x_0}M$ with the origin in $x_0$ and the tips on the indicatrix segment $N_{\tau_2}^{\tau_1}:=\{N(\tau):\tau\in[\tau_1,\tau_2]\}$ and to join the 
points $u_0^-:=(x_0,N^-)$ and  $u_0^+:=(x_0,N^+)$ in $\Sigma$ by the arc of indicatrix curve $(x_0,N_{\tau_2}^{\tau_1})$. Unfortunately,  
this method is not yet good enough because one does not obtain in this way a smooth section of $TM$ along $\gamma$ and therefore, the existence of the prolongation vector field $X$ is not guaranteed anymore.

However, this idea works well if we consider a smooth variation of $\gamma$ on $M$.

Indeed, let us consider a variation $\tilde{\gamma}_\varepsilon:[0,1]\to M$ of $\gamma$ depending on a small $\varepsilon>0$ such that $\lim _{\varepsilon \to 0} \tilde{\gamma}_\varepsilon =\gamma$ as set of points. 

We define 
\begin{equation}\tag{5.5}\label{5.5}
\tilde{\gamma}_1(\varepsilon,t)=\exp_{\gamma(t)}(\varepsilon N(t)),
\end{equation}
where $N$ is the normal vector field along $\gamma$. Since $N$ has a discontinous jump from $N^-$ to $N^+$ at $x_0$, the curve $\tilde{\gamma}_1$ will also have a jump.

Indeed, let us remark that for a fixed small enough $\varepsilon >0$, we obtain a smooth curve  $\tilde{\gamma}_1(\varepsilon,t)$ on $M$ going around $\gamma$, while for a fixed $t$ we have a transversal Finslerian geodesic with initial conditions $(\gamma(t),N(t))$.

Remark also that for a fixed $t_1$ the tangent vector $\tilde{T}_{t_1}(\varepsilon)$ of 
  $\tilde{\gamma}_1$  at the point  $\tilde{\gamma}_1(t_1,\varepsilon)$ 
is given by the parallel translation of the tangent vector $T(t_1)$ 
of  $\gamma$ at the point $\gamma(t_1)$ along the transversal geodesic 
$\exp_{\gamma(t_1)}(\varepsilon N(t_1))$, where $g_{N(t_1)}(N(t_1),T(t_1))=0$. Using now the properties of parallel displacement (see \cite{BCS2000} p. 140, 141) it follows that at any small enough $\varepsilon>0$ we have
\begin{equation*}
g_{\tilde N_{t_1}(\varepsilon)}
(\tilde N_{t_1}(\varepsilon),\tilde T_{t_1}(\varepsilon))
=0,
\end{equation*}  
where $\tilde N_{t_1}(\varepsilon)$ is the tangent vector of the transversal geodesic 
$\exp_{\gamma(t_1)}(\varepsilon N(t_1))$ at the point $\gamma_1(t_1,\varepsilon)$

These remarks assure us that the variation curve  $\tilde{\gamma}_1$ has its ends on the geodesics $\exp_{x_0}(\varepsilon N^-)$ and $\exp_{x_0}(\varepsilon N^+)$, for small enough $\varepsilon >0$, i.e.  $\tilde{\gamma}_1$ is not a closed loop.

Next, we will complete the curve  $\tilde{\gamma}_1$ with an arc of curve 
$\tilde{\gamma}_2$ that connects smoothly the ends of   $\tilde{\gamma}_1$ such that  $\tilde{\gamma}= \tilde{\gamma}_1\cup  \tilde{\gamma}_2$ is a closed smooth variation of $\gamma$ on $M$. The easiest way to do this is exponentiate the indicatrix arc between $N^-$ and $N^+$, i.e. we consider 
\begin{equation}\tag{5.6}\label{5.6}
 \tilde{\gamma}_2(\varepsilon)=\exp_{x_0}(\varepsilon N_{\tau_2}^{\tau_1}).
\end{equation} 

One can now easily see that $\tilde{\gamma}= \tilde{\gamma}_1\cup  \tilde{\gamma}_2$ is a closed smooth variation near $\gamma$ whose tangent vector
$\tilde{T}$ is given along $\tilde{\gamma_1}$  by the parallel displacement of $T$ along the transversal geodesic $\sigma_t(\varepsilon)={\tilde{\gamma}_1}(\varepsilon,t)$, and along $\tilde{\gamma_2}$ by $\exp_{x_0*}W$, where $W$ is the tangent vector along the indicatrix curve. Gauss Lemma for Riemann--Finsler manifolds (see for example \cite{BCS2000}, p. 140) assures us that 
$\hat{g}_{x_0}(\varepsilon N,W)=0$ and $g_{\tilde N}(\tilde N,\tilde{T_2})=0$, where $\tilde N$ and $\tilde T_2$ are the normal and tangent vectors along $\tilde \gamma_2$, respectively. 

From the discussion above, one can see now that the tangent vector of 
$\tilde \gamma_1$ at 
$x_0^-=\exp_{x_0}(\varepsilon N^-)$ is $g_{\tilde N^-}$ orthogonal to 
$\tilde N^-:=\frac{d}{d\varepsilon}\tilde \gamma_1(t,\varepsilon)$ and that the tangent vector of  $\tilde \gamma_2$ at the same point  $x_0^-$ is also $g_{\tilde N^-}$ orthogonal to 
$\tilde N^-$ due to Gauss Lemma, therefore the unitary left and right tangent vectors at $x_0^-$ have the same direction, so they must coincide (see Figure 4).

Therefore we can conclude that the curve $\tilde \gamma$ is smooth at $x_0^-$ when we take the limit $\varepsilon \to 0$. The same argument applies at $x_0^+=\exp_{x_0}(\varepsilon N^+)$.

We point out however that since we have moved the point $x_0$ a little along the transversal geodesic $\exp_{x_0}(\varepsilon N^-)$ the indicatrix also changes from $\Sigma_{x_0}$ to $\Sigma_{x_0^-}$. However, we will finally take the limit  $\varepsilon \to 0$ so this small displacement cannot cause much harm.

\smallskip
\begin{center}
\includegraphics[height=6cm, angle=0]{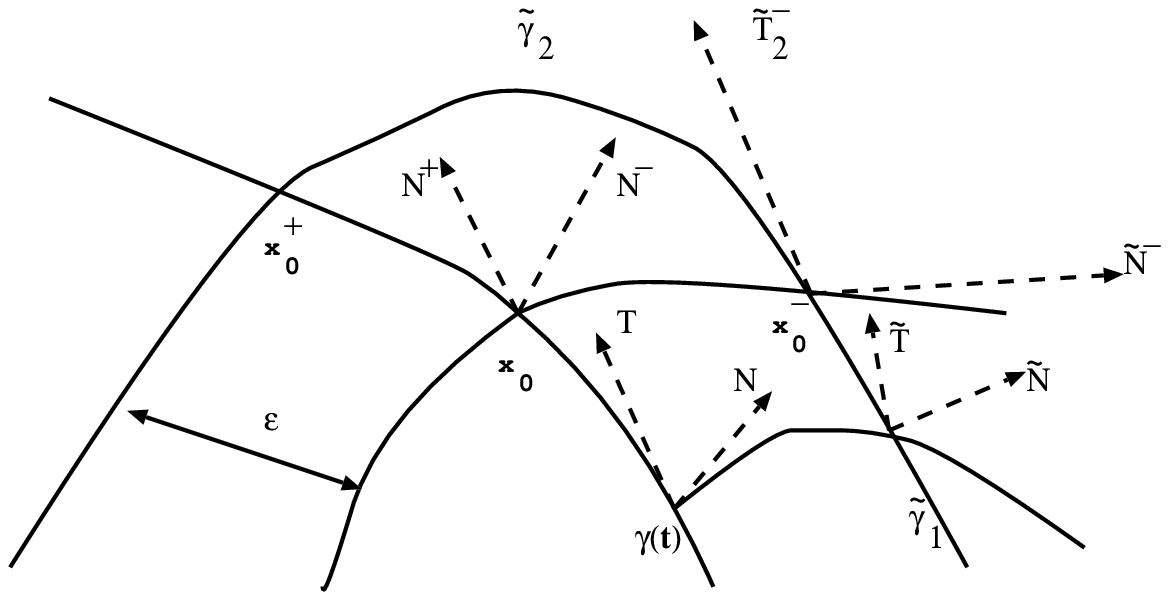}
\end{center}
{\bf Figure 4.} A magnified view of the landscape around the point $x_0$. 
\smallskip

Having now all these done, we can now consider the bounded domain $\tilde{D}\subset M$ with smooth boundary $\partial \tilde{D}=\tilde{\gamma}=\tilde{\gamma}_1\cup\tilde{\gamma}_2$ and apply to it the same method as in \S4.

Indeed, writing our Topological Lemma for $\tilde{D}$ and taking the limit, we obtain
\begin{equation*}
 - \int_{D}{X}^*(d\Pi)+\lim_{\varepsilon\to 0} \int_{\tilde{N}(\partial {\tilde D})}\Pi= \mathcal X( D),
\end{equation*}
with the obvious notations. 

The term concerning the boundary becomes
 \begin{equation*}
\begin{split}
\lim_{\varepsilon\to 0} \int_{\tilde{N}(\partial {\tilde D})}\Pi & =
\lim_{\varepsilon\to 0} \int_{\tilde{\gamma_1}}\tilde{N}^*(\Pi)+
\lim_{\varepsilon\to 0} \int_{\tilde{\gamma_2}}\tilde{N}^*(\Pi)\\
& =
\int_{\gamma}{N}^*(\Pi)+\frac{1}{L(x_0)}\int_{N_{\tau_2}^{\tau_1}}\omega_1^{\ 2}.
\end{split}
\end{equation*}

We are going to compute now the second integral in the sum above.

Remark that we are now integrating on the segment $N_{\tau_2}^{\tau_1}$
where there is no variation of $x$, therefore the integrand reads
 \begin{equation*}
\omega_1^{\ 2}=\frac{\sqrt{g}}{F^2}(y^1\delta y^2-y^2\delta y^1)=
\frac{\sqrt{g}}{F^2}(y^1 \dot{y}^2-y^2 \dot{y}^1)d\tau,
\end{equation*}
where $y^i=y^i(\tau)$ is a unit speed parameterization of the indicatrix 
$\Sigma_{x_0}$, and $\{dx^i,\frac{1}{F}\delta y^i\}$ is the dual cobasis of 
the adapted basis $\{\frac{\delta}{\delta x^i},F\frac{\partial}{\partial y^i}\}$.
Here $\delta y^i=dy^i+N^i_jdx^j$ (see \cite{BCS2000}, p. 96 for details).

Recall that the tangent vector to the indicatrix is given
by
\begin{equation*}
\lambda=\dot{y}^1\frac{\partial}{\partial y^1}+
\dot{y}^2\frac{\partial}{\partial y^2}.
\end{equation*}

Therefore, we have
\begin{equation}\tag{5.7}\label{5.7}
\begin{split}
\int_{N_{\tau_2}^{\tau_1}}\omega_1^{\ 2} & =\int_{\tau_1}^{\tau_2}\omega_1^{\
 2}(\lambda)d\tau=\int_{\tau_1}^{\tau_2}\frac{\sqrt{g}}{F^2}(y^1
 \dot{y}^2-y^2 \dot{y}^1)d\tau \\
& =\int_{\tau_1}^{\tau_2}d\theta=\measuredangle_{x_0} (N^-,N^+).
\end{split}
\end{equation}

Putting all these together, we obtain the following main result

\bigskip

{\bf Theorem 5.1.} (The Gauss--Bonnet theorem for Finsler surfaces with regular
$C^\infty$ piecewise boundary)

{\it Let $(M,F)$ be a compact oriented Finslerian surface and $D\in M$ a domain with regular piecewise $C^\infty$ boundary $\partial D=\gamma$, that consists of the union of $k$ piecewise smooth curves. Let $N:\partial D\to \Sigma$ be the inward pointing Finslerian unit normal on $\partial D$. 

Then, we have
\begin{equation}\tag{5.8}\label{5.8}
\begin{split}
 \int_D \cfrac{1}{L(x)} & \Bigl[X^*(K\ \omega^1\wedge \omega^2
  - J\omega^1\wedge\omega^3)-d\log L(x)\wedge X^*(\omega^3)
\Bigr]\\
& +\int_{\gamma}\frac{1}{L(x)}\frac{k_T^{(N)}(t)}{\sigma(t)}dt
+\sum_{s=1}^k\frac{1}{L(x_s)}\measuredangle_{x_s} (N^-,N^+)=
\mathcal X(D),
\end{split}
\end{equation}
where 
$L(x)$ is the Riemannian length of the indicatrix $\Sigma_x$,
$X$ is a unit prolongation of $N$, 
$K$ is the Gauss curvature, 
$\measuredangle_{x_s} (N^-,N^+)$ the Landsberg angle of the unit vectors $N^-$ and $N^+$,
and $\mathcal X(D)$ is the Euler characteristic of $D$. 
}

\bigskip

{\bf Remarks.}
\begin{enumerate}
\item
If $(M,F)$ is a Riemannian manifold, then the Gauss--Bonnet theorem
formulated above reduces to the classical Gauss-Bonnet theorem on
Riemannian manifolds. Indeed, it suffices to remark that, in the
Riemannian case, the Euclidean angle $\measuredangle_{x_s} (T^-,T^+)$
equals the angle $\measuredangle_{x_s} (N^-,N^+)$ which is also an
Euclidean angle. Nevertheless, in the Riemannian case, the sum of
interior and exterior angles at a corner equals $\pi$, but this is not
the case anymore in the Finslerian case as already discussed in $\S2$.

\item If $D$ is a domain with regular piecewise $C^\infty$ boundary $\partial D=\gamma$ on a Landsberg surface $(M,F)$, then we obtain
\begin{equation}\tag{5.9}\label{5.9}
 \cfrac{1}{L}\int_D K\sqrt{g} dx^1\wedge dx^2
 +\cfrac{1}{L}\int_{\gamma}\frac{k_T^{(N)}(t)}{\sigma(t)}dt+
\frac{1}{L} \sum_{s=1}^k\measuredangle_{x_s} (N^-,N^+)
=
\mathcal X(D),
\end{equation}
with the obvoius notations.

\end{enumerate}



\section{A Hadamard type theorem for $N$-parallels}

We are going to discuss here an application of the Gauss-Bonnet formula (5.9) for Landsberg surfaces.

In Riemannian geometry it is known that the Gauss-Bonnet theorem imposes restrictions on the behavior of geodesics. Namely, Hadamard theorem states that on a simply connected Riemannian surface of nonpositive Gauss curvature $K\leq 0$, a geodesic cannot have self intersections.

We are going to prove a similar result for the $N$-parallels of a Landsberg surface. First, remark the following

{\bf Lemma 6.1.}

{\it Let $x_0$ be a point on M, and let us denote by $\Sigma_{x_0}\in T_{x_0}M$  the indicatrix curve of $(M,F)$ at $x_0$. Then we have
\begin{equation}\tag{6.1}\label{6.1}
 \frac{1}{L(x_0)}\ \measuredangle_{x_s} (N^-,N^+)<1,
\end{equation}
where $L(x_0)$ is the Riemannian length of the indicatrix $\Sigma_{x_0}$ and $\measuredangle_{x_s} (N^-,N^+)$ is the Landsberg angle of the unit length vectors $N^-$, $N^+$.
}

The proof is trivial. For a positive orientation, the Riemannian length of the indicatrix arc ${{\stackrel{\frown}{N^-N^+}}}$ at $x_0$ is always smaller than the total length of the indicatrix  $\Sigma_{x_0}$ (see Figure 3).

We can give now

{\bf Theorem 6.2.} 

{\it 
On a simply connected Landsberg surface $(M,F)$ of nonpositive Gauss curvature $K\leq 0$, the N-parallels cannot have self-intersections.
}

{\it Proof.} Let us assume that the $N$-parallel $\gamma:[a,b]\to M$ can have self intersections, and let us denote such a point by $x_0$. 

This is equivalent with saying that on $M$ we have a domain $D$ with close regular piecewise $C^\infty$
boundary $\partial D=\gamma$. The boundary curve on $M$ is an $N$-parallel having a corner at $x_0$.

\smallskip
\begin{center}
\includegraphics[height=5cm, angle=0]{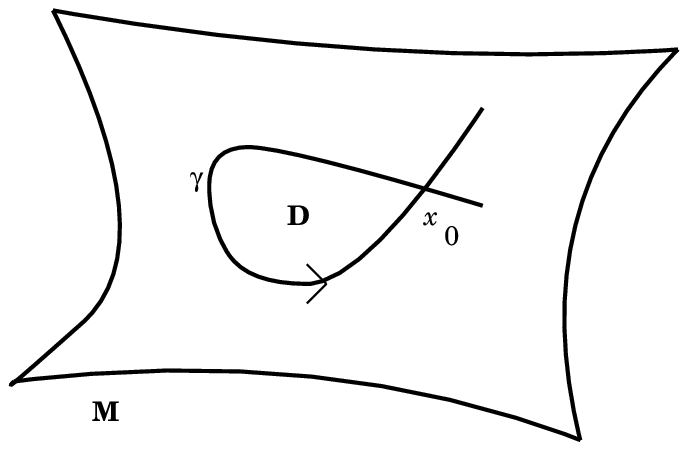}
\end{center}
{\bf Figure 5.} A self intersecting $N$-parallel curve.
\smallskip

Applying now the Gauss-Bonnet formula  (5.9) for the domain $D$ with boundary $\partial D=\gamma$ we obtain
 \begin{equation}\tag{6.2}\label{6.2}
\cfrac{1}{L}\int_D K\sqrt{g} dx^1\wedge dx^2+
\frac{1}{L} \measuredangle_{x_0} (N^-,N^+)
=1,
\end{equation}
where $N^-$, $N^+$ are the left and right normals, respectively, to the boundary in the 
point $x_0$ as before.

One can see now that this formula leads to a contradiction showing in this way that the assumption is false. Indeed, since $K\leq 0$ is a nonpositive function, the integral in the left hand side of (\ref{6.2}) is nonpositive. On the other hand, from Lemma 6.1. we know that 
the second term in the sum in the left hand side of (\ref{6.2}) is less than 1. But this is not possible, therefore we have reached to a contradiction.

It follows that the $N$-parallel curve $\gamma$ cannot have a self intersection, in other words, the situation on Figure 5 cannot happen.

\begin{flushright}
 Q. E. D. 
\end{flushright}

\bigskip

{\bf Remarks.}
\begin{enumerate}
 \item Recall that {\it Euler's theorem for polyhedra} states that for any triangulation of a compact surface $M$, the Euler characteristic of $M$ is given by
\begin{equation}
 \mathcal X(M)=\sharp vertices - \sharp edges + \sharp faces,
\end{equation}
where the symbol $\sharp$ ``means the number of''. 
In particular, if we have a bounded region $D$ on a simply connected surface $M$ like in Figure 5, then $D$ is homeomorphic to a triangle, i.e.
\begin{equation}
 \mathcal X(D)=\sharp vertices - \sharp edges + \sharp faces=3-3+1=1.
\end{equation}
This is the reason we have 1 in the right hand side of (\ref{6.2}).
\item There is a second part of the Hadamard theorem that states that on any simply connected Riemannian surface of non-positive Gauss curvature $K\leq 0$ two distinct geodesics cannot have two points of intersection. This kind of result also extends to the case of $N$-parallels, but it is a little more complicated and is going to be discussed in a forthcoming paper together with other applications of the Gauss-Bonnet theorem.
\end{enumerate}



\section{Appendix: The existence and unicity of $N$-parallels}

$\quad$ Formula (\ref{3.16}) is useful for the study of existence and unicity of the $N$-parallels of a Finsler surface $(M, F)$.\\

Indeed, following an idea of M. Matsumoto \cite{M1984} from the conditions that define the $N$-parallels, namely 
\begin{equation*}
F(x, N) = 1, \qquad g_N(N, T) = 0,
\end{equation*}
or, equivalently,
\begin{eqnarray*}
& &g_{ij}(x, N)\cdot N^iN^j = 1 , \\
& &g_{ij}(x, N)\cdot N^iT^j = 0,
\end{eqnarray*}
where $i, j = 1, 2$, it follows 
\begin{eqnarray*}
\big\lbrack g_{i1}(x, N)\cdot N^i\big\rbrack T^1 + \big\lbrack g_{i2}(x, N)\cdot N^i\big\rbrack T^2 = 0.
\end{eqnarray*}

From here, it follows that there exists a positive scalar $k$ such that 
\begin{eqnarray*}
-\frac{ g_{i1}(x, N)\cdot N^i}{T^2} = \frac{ g_{i2}(x, N)\cdot N^i}{T^1} = k > 0
\end{eqnarray*}
(or with opposite signs) and therefore,
\begin{eqnarray*}
& & g_{i1}(x, N)\cdot N^i = -k\cdot T^2 , \\
& & g_{i2}(x, N)\cdot N^i = k\cdot T^1.
\end{eqnarray*}

Using now the $0$-homogeneity of $g_{ij}$, we obtain the equations 
\begin{equation}\tag{7.1}\label{7.1}
 \begin{split}
&g_{i1}(x, p)\cdot p^i = - T^2 , \\
&g_{i2}(x, p)\cdot p^i =  T^1,
 \end{split}
\end{equation}
where $i = 1, 2$ and $p$ is a vector proportional to $N$.\\

Taking into account that Jacobian of the equation $\eqref{7.1}$ is just 
\begin{equation*}
det\vert\, g_{ij}(x,p)\vert \not= 0
\end{equation*}
it follows by the Theorem of Implicit Functions that we can solve these equations with respect to the unknowns $p^1$, $p^2$.\\

Finally, we can put 
\begin{equation}\tag{7.2}\label{7.2}
N^i := \frac{p^i}{F(x,p)}, \qquad i = 1, 2.
\end{equation}

One can easily see that this $N = (N^i)$ satisfies condition $\eqref{3.1}$.\\

We point out that the solutions $N^1$, $N^2$ of the equation $\eqref{7.1}$ depend actually on $T$.\\

We can rewrite $\eqref{3.16}$ as 
\begin{equation}\tag{7.3}\label{7.3}
\frac{d^2\gamma^i}{dt^2} + \Gamma^i_{jk}\big(\gamma(t), N(t)\big)\frac{d\gamma^j}{dt}\frac{d\gamma^k}{dt} = 
\frac{d}{dt}\Bigl[ \log \sigma(t) \Bigr] \frac{d\gamma^i}{dt},
\end{equation}
where $N(t) = N\big(\gamma(t), \dot\gamma(t)\big)$ from $\eqref{7.1}$.\\

An initial condition can be given by
\begin{equation}\tag{7.4}\label{7.4}
\begin{split}
\gamma^i(t_0) = 0,\\
\dot\gamma^i(t_0) = T_0^i,
\end{split}
\end{equation}
with $i = 1, 2$ and corresponding the normal initial condition
\begin{equation}\tag{7.5}\label{7.5}
\begin{split}
\gamma^i(t_0) = 0,\\
N^i(t_0) = N_0,
\end{split}
\end{equation}
where $N_0$ are given as solutions of $\eqref{7.1}$ for $T = T_0$.\\

Then, by a similar argument as in the case of geodesics, we know from the general theory of ODEs that $\eqref{7.3}$ with initial conditions $\eqref{7.4}$ have unique solutions.

A detailed study of the $N$-parallels will be given elsewhere.

\bigskip

{\bf Acknowledgments.}$\;$ 
The authors would like to express their gratitude to Professors L. Tamassy and Gh. Pitis for many useful discussion on this topic. We thank to Professor K. Shiohama for bringing this topic into our attention and his continuous encouragement. Finally, we are grateful to the referee for the careful reading of the paper and for many helpful suggestions.


\bigskip

\begin{flushleft}
Jin-Ichi Itoh\\
Kumamoto University\\
Kumamoto, Japan\\
E-mail: j-itoh @ kumamoto-u.ac.jp
\end{flushleft}

\begin{flushleft}
Sorin V. Sabau\\
Hideo Shimada\\
Tokai University, Sapporo Campus\\
Sapporo, Japan\\
E-mail: sorin @ tspirit.tokai-u.jp\\
E-mail: shimadah @ tokai-u.jp
\end{flushleft}


\begin{thebibliography}{[BCS2000]}

\bibitem[B2007]{B2007}
     {Bao},~D., 
{\it On two curvature-driven problems in Riemann-Finsler geometry}, 
Finsler geometry, Sapporo 2005, Advances Studies in Pure Math., {\bf 48} (2007), 
19--71.

\bibitem[BC1996]{BC1996}
     {Bao},~D., {Chern}, ~S.S.,
{\it A Note on Gauss-Bonnet Theorem for Finsler spaces}, Annals of Math., 
{\bf 143}(1996), 233--252. 

\bibitem[BCS2000]{BCS2000}
     {Bao},~D., {Chern}, ~S.S., {Shen}, ~Z.,
     {\it An Introduction to Riemann Finsler Geometry}, Springer, GTM 200, 
2000.

\bibitem[BS1994]{BS1994}
     {Bao},~D., {Shen}, ~Z.,
{\it On the volume of unit tangent spheres in a Finsler manifold}, Results 
in Math., {\bf 26}(1994), 1--17.

\bibitem[Br1997]{Br1997}
	       {Bryant},~R.,
{\it Projectively flat Finsler 2-spheres of constant curvature}, Selecta
		 Math. (N.S.), {\bf vol. 3, no. 2} (1997), 161--203.

\bibitem[Br2002]{Br2002}
	       {Bryant},~R.,
{\it Some remarks on Finsler manifolds with constant flag curvature},
Houston Journal of Mathematics, {\bf vol. 28, no.2} (2002), 221--262.  

\bibitem[I1978]{I1978}
     {Ichijyo},~Y.,
{\it On special Finsler connections with vanishing $hv$-curvature tensor}, 
Tensor N.S., {\bf 32}(1978), 146--155.

\bibitem[L2002]{L2002}
	{Lackey},~B.,
{\it On the Gauss-Bonnet formula in Riemann-Finsler geometry}, Bull. London Math. Soc. {\bf 34} (2002), 329--340.

\bibitem[Ma2008]{Ma2008}
{Matveev},~S.V.,
{\it On "Regular Landsberg metrics are always Berwald" by Z.I.Szabo},
arXiv:0809.1581v1 $\lbrack$ math.DG$\rbrack$ 9 Sept.2008.

\bibitem[M1984]{M1984}
     {Matsumoto},~M.,
     {\it Theory of $Y$-extremal and minimal hypersurfaces in a Finsler space. On Wegener's and Barthel's theories}, 
J. Math. Kyoto Univ., {\bf 26(4)} (1986), 647-665,

\bibitem[M1986]{M1986}
     {Matsumoto},~M.,
     {\it Foundations of Finsler Geometry and Special Finsler Spaces}, 
Kaiseisha Press, Otsu, Japan, 1986.

\bibitem[Mil1965]{Mil1965}
     {Milnor},~J.,
     {\it Topology from the differentiable viewpoint},
University Press of Virginia, Charlotesville, 1965.


\bibitem[SS2007]{SS2007}
     {Sabau},~V.,~S., {Shimada},~H.,
{\it Riemann-Finsler surfaces}, 
Finsler geometry, Sapporo 2005, Advances Studies in Pure Math., {\bf 48} (2007), 
125--162.

\bibitem[Sh1996]{Sh1996}
	{Shen},~Z.,
{\it Some formulas of Gauss-Bonnet-Chern type in Riemann-Finsler geometry},
J. reine angew. Math., {\bf 475} (1996), 149--165.

\bibitem[Sh2001]{Sh2001}
	{Shen},~Z.,
{\it Lectures on Finsler Geometry},
World Scientific, 2001.

\bibitem[SST2003]{SST2003}
	{Shiohama}, ~K., {Shioya}, ~T., {Tanaka}, ~M.,
{\it The geometry of total curvature on complete open surfaces},
Cambridge University Press, 2003.

\bibitem[Spiv1979]{Spiv1979}
         {Spivak}, ~M.,
{\it A Comprehensive Introduction to Differential Geometry, Second Edition},
{\bf Vol. V},
Publish or Perish, Inc., 1979.

\bibitem[Sz1981]{Sz1981}
     {Szab\'{o}},~Z.,
{\it Positive definite Berwald spaces (Structure theorems on Berwald 
spaces)}, Tensor N.S., {\bf 35}(1981), 25--39.

\bibitem[Sz2008a]{Sz2008a}
{Szabo}, ~Z.,
{\it All regular Landsberg metrics are Berwald}, to appear in AGAG (2008), on line publication: 
http:\slash \slash dx.doi.org \slash 10.1007\slash s10455-008-9115-y.


\bibitem[Sz2008b]{Sz2008b}
{Szabo}, ~Z.,
{\it Correction to ``All regular Landsberg metrics are Berwald''}, preprint, 2008.
\end{thebibliography}
\end{document}